\documentclass{birkjour}
 \newtheorem{thm}{Theorem}[section]
 
 \newtheorem{lem}[thm]{Lemma}
 \newtheorem{prop}[thm]{Proposition}
 \theoremstyle{definition}
 \newtheorem{defn}[thm]{Definition}
 \theoremstyle{remark}
 \newtheorem{rem}[thm]{Remark}
 
 \numberwithin{equation}{section}

\begin{document}

%
%
%
%
%
%
%
%
%

\title[Entire solutions]
 {Asymptotic behavior of entire solutions for degenerate partial differential inequalities on Carnot-Carath\'{e}odory metric spaces and \\ Liouville type results}

\author[Markasheva]{Markasheva V.A.}

\address{%
Universit\'{a} di Bologna\\
Dipartimento di Matematica,\\
Piazza di Porta San Donato, 5\\
BO 40127 Bologna \\
Italy}

\email{vira.markasheva@unibo.it, w9071981@gmail.com}

\thanks{Research supported by the program
"Erasmus Mundus Action 2 Lot 5 MID" provided by the European
Commission and coordinated by the University of Turku (Finland),
grant No. MID2012B207. The author has been additionally funded by
the project A.G.A.P.E. (Analysis in lie Groups and Applications to
Perceptual Emergences) coordinated by the Department of Mathematics
of the University of Bologna.}



\date{09.05.2015}

\begin{abstract}
 This article is devoted to the study of the behavior of
generalized entire solutions for a wide class of quasilinear
degenerate inequalities modeled on the following prototype with
p-Laplacian in the main part
\begin{equation*}
{\underset{m}{\overset{i=1}{\sum}}} X_i^*(|\mathbf{X}u|^{p-2} X_i
u)\geq |u|^{q-2}u, \ \ x\in {\mathbb{R}}^{n},\ q>1,\ p>1,
\end{equation*}
where ${\mathbb{R}}^{n}$ is a Carnot-Carath\'{e}odory metric space,
generated by the system of vector fields
$\mathbf{X}=(X_1,X_2,..,X_m)$ and $X_i^*$ denotes the adjoint of
$X_i$ with respect to Lebesgue measure. For the case where $p$ is
less than the homogeneous dimension $Q$ we have obtained a sharp a
priori estimate for essential supremum of generalized solutions from
below which imply some Liouville-type results.
\end{abstract}

\maketitle
\section{Introduction}
The research subject analyzed in this paper is the asymptotic
behavior of non-negative entire solutions of the following
subelliptic differential inequality

\begin{equation}\label{mark4.1}
{\overset{m}{\underset{i=1}{\sum}}} X_i^* A_i(x,u,\mathbf{X} u) \geq
f(x,u),  \ \ x\in R^n,
\end{equation}
 under the conditions:
$$
\bf{(A_1)}\begin{array}{l} A(x,u, \xi)\  \mbox{is a Carnot
function},\\\mbox{such that for all} \ \xi\in {\mathbb{R}}^{m}:
\nu_1 |{\xi}|^p \leq A (x,u,{\xi}){\xi}\  \mbox{on}\
{\mathbb{R}}^{n},
 \end{array}
$$
$$
\bf{(A_2)}\begin{array}{l} A(x,u,\xi)\ \mbox{is a Carnot
function},\\\mbox{such that for all}\ \xi\in {\mathbb{R}}^{m}: A
(x,u,{\xi})\leq \nu_2 |{\xi}|^{p-1}\  \mbox{on}\ {\mathbb{R}}^{n},
 \end{array}
$$
$$
\bf{(A_3)}\begin{array}{l} f(x,u)\ \mbox{is a Carnot
function},\\\mbox{such that for all}\ u \in {\mathbb{R}}: f(x,u)
\geq \nu_3|u|^{q-2}u\  \mbox{on}\  {\mathbb{R}}^{n},
 \end{array}
$$
where $q>1$, $\ 1<p<Q$, ${\mathbb{R}}^{n}$ is a
Carnot-Carath\'{e}odory metric space of homogeneous dimension $Q$,
generated by the system of vector fields
$\mathbf{X}=(X_1,X_2,..,X_m)$ and $X_i^*$ denotes the adjoint of
$X_i$ with respect to Lebesgue measure. $d_{CC}(x)$ is
Carnot-Carath\'{e}odory metric distance between $x$ and $0$.  As it
was mentioned in the abstract, a typical example of \eqref{mark4.1}
is the following
\begin{equation*}
{\underset{m}{\overset{i=1}{\sum}}} X_i^*(|\mathbf{X}u|^{p-2} X_i
u)\geq |u|^{q-2}u, \ \ x\in {\mathbb{R}}^{n},\ q>1,\ 1<p<Q.
\end{equation*}
Both inequalities are considered under the restriction $u\geq 0$. To
describe what kind of results we can expect from such type elliptic
problems, it would be interesting to formulate in Euclidean settings
some well-known facts. Let us consider the following equation
\begin{equation}\label{mark4.140}
 div(|\nabla u|^{p-2}\nabla u )= \nu |u|^{q-2}u, \ \ x\in {\mathbb{R}}^{n},\ 1<q,\
1<p<n.
\end{equation} Here ${\mathbb{R}}^{n}$ is Euclidean space and $\nabla u$ means a classical notion of gradient of a function $u$.
The simple consideration below illustrate the fact that the
asymptotic behavior of entire solutions of this equation tightly
depends on exponents $p$ and $q$.

Let $q<p$. There exists a positive constant $\nu$ such that
$u(x)=|x|^{p/(p-q)}$ is an entire solution of the equation
\eqref{mark4.140} and $\triangle_p u= \nu u^{q-1}$ a.e.(the last
fact explains the assumption $q>1$). This trivial example becomes
more interesting in comparison with some Liouville type results
obtained by Serrin in \cite{Ser2009} for $C^1$ smooth solutions: if
$1<q<p$ and
$u(x)=\bar{\bar{o}}(|x|^{p/(p-q)})$ 
 then $u\equiv 0$ (the
description of these results is done according to the notations of
this paper and to the case). The example shows that Serrin's
condition is sharp and cannot be improved. The closer to $p$ we
choose an exponent $q$ in the equation the faster it's solutions
must grow at infinity.

In case $q=p$ we could expect that solutions grow faster than any
algebraic function and this fact was proved by Serrin in the same
paper. Indeed, as examples of solutions, one can consider functions
$\exp{(x_i)},\ i=1,2,..,n.$ The method used by Serrin is built on
earlier ideas of Redheffer \cite{Red1960} in combination with
Serrin's ideas. The approach supposes  $C^1$ continuity of solutions
as an important ingredient. A curious reader can find an additional
information about some preceding investigations in Euclidean
settings also in papers of Farina \cite{Far2007}, Tkachev
\cite{Tkach1994}, Brezis \cite{Brez1984}, Benguria, Lorca and Yarur
\cite{BenLorYar1994}, Haito, Usami \cite{NaitUsam1997}.

For the case $q>p$ Serrin proved that there could be only trivial
solution $u\equiv 0.$ This result but for inequalities was treated
also by Mitidieri, D'Ambrosio \cite{D'AmbrMit2010} as an auxiliary
result under Euclidean settings and by D'Ambrosio \cite{D'Ambr2009}
on Lie groups.

The method we are going to apply in case $q\leq p$ to prove
Liouville type result  was introduced in the article of Kondrat'ev
and Landis \cite{KondrLand1988} for the case $q>p$ in 1988 and was
carefully developed and generalized by Mitidieri and Pohozhaev  for
the another type of inequalities $-L u\geq f(x,u)$ in their seminal
article \cite{MitPoh2001}.

Moreover, the apriori estimate of the maximum of solutions from
below is non-trivial. This idea is new even in Euclidean settings.
It belongs to Analoti Tedeev, who has obtained the above-mentioned
estimate for non-negative solutions of Cauchy problem for parabolic
equations with p-Laplacian on Grushin metric spaces (see
\cite{Ted1}). Solutions of Cauchy problem for parabolic equations on
the space are  distributing their mass all over the space with time.
The maximum of non-negative solutions is tending to zero. The bound
from below in this situation proves that the estimate of maximum is
precise.

By the contrary, in the present case we have no decay. We are
dealing with entire solutions which are growing close to infinity.
The estimate of the maximum from below is much more important than
from above.

Let us formulate our main results:

\begin{thm}\label{markt4.1}
Let $1 <q<p<Q$ and let $u$ be a weak non-negative entire solution of
the differential inequality \eqref{mark4.1} then
\begin{itemize}
\item[i)] $\exists \  C_1>0\ $ and $\ R_0>0\ $ such that for $\ \forall\  R\geq
R_0$ on $d_{CC}-$annulus $A_{R/2}^{2R}$
\begin{equation}\label{mark4.3}
||u||_{\infty, A_{R/2}^{2R}}\geq C_1 R^{\frac{p}{p-q}},
\end{equation}
\item[ii)] if $\ u={\bar{o}}(d_{CC}^{p/(p-q)}(x))$ as $d_{CC}(x)\rightarrow
+\infty\ $ then $u\equiv 0.$
\end{itemize}
\end{thm}
\begin{rem}\label{markr4.1}
When the a priori estimate \eqref{mark4.3} is proved it is evident
that the second assertion of the theorem follows immediately.
\end{rem}
\begin{thm}\label{markt4.2}
Let $1 <q=p<Q$ and let $u$ be a weak non-negative entire solution of
the differential inequality \eqref{mark4.1} then
\begin{itemize}
\item[i)] $\exists \  C_2>0\ $ and $\ R_0>0\ $ such that for $\ \forall\  R\geq
R_0$
\begin{equation}\label{mark4.4}
||u||_{\infty, A_{R/2}^{2R}}\geq C_2 R^{\frac{p}{p-1}},
\end{equation}
\item[ii)] Let $\gamma:\ 0<\gamma \leq p/(p-1)$ be an any fixed real number then if  $\ u={\bar{o}}(d_{CC}^{\gamma}(x))$ as $d_{CC}(x)\rightarrow
+\infty\ $ then $u\equiv 0.$
\end{itemize}
\end{thm}
\begin{rem}\label{markr4.1} The example of Serrin shows us that the
estimate \eqref{mark4.4} is not precise. ii) must be true for any
positive $\gamma$. We are going to prove the assertion \textit{ii)}
for the $\gamma \geq p/(p-1)$ in our next paper.
\end{rem}
\begin{thm}\label{markt4.3}
Let $1\leq p<q$ and let $u$ be a weak entire solution of the
differential inequality \eqref{mark4.1} then $u\equiv 0.$
\end{thm}
\begin{rem}\label{markr4.3}
We added this theorem into the paper for the completeness of the
picture with the exponents, from one point of view, and because the
proof of this result is unexpectedly simple and short using the
Kondrat'iev-Landis-Mitidieri-Pohozhaev method, from another point of
view.
\end{rem}
\begin{rem}\label{markr4.4}
It is also worth to mention that we can formulate  Theorem
\ref{markt4.3} not only for p-Laplacian case ($1\leq p<\infty$) but
for all differential operators with non-negative characteristic
forms, for example, the theorem is true for  mean curvature operator
and  for total variation operator as well.
\end{rem}

The structure of the paper is as follows: in the second section we
 will present some auxiliary statements  and  introduce the notion of a generalized entire solution to the
 equation \eqref{mark4.1}. The  third section is devoted to the proof of statements \textit{i)}
 both in Theorem \ref{markt4.1} and in Theorem \ref{markt4.2}.
 In the fourth section, we prove Theorem \ref{markt4.3}.

\bigskip
\section{Auxiliary results}

We consider in ${\mathbb{R}}^{n}$ a system
$\mathbf{X}=(X_1,X_2,..,X_m)$ of vector fields
\begin{equation*}
X_j= \underset{k=1}{\overset{n}\sum}
b_{jk}(x)\frac{\partial}{\partial x_k},\ \ j=1,..,m,
\end{equation*}
having  real-valued, locally Lipschitz-continuous $b_{jk}$. Through
the paper if $u$ is a non smooth  then $X_j u$ will be meant in the
distributional sense. For the system $\mathbf{X}=(X_1,X_2,..,X_m)$
we denote by $|\mathbf{X} u|= \left(\underset{j=1}{\overset{m}\sum}
(X_j u)^2\right)^{1/2}$ the length of the horizontal gradient
$$\mathbf{X} u =(X_1 u,X_2 u,..,X_m u).$$
Thus let us consider also ${\mathbb{R}}^{n}$ as a
Carnot-Carath\'{e}odory metric space of a homogeneous dimension $Q$
with Carnot-Carath\'{e}odory metric distance $d_{CC}(x)$ generated
by the system of vector fields $\mathbf{X}=(X_1,X_2,..,X_m)$ as it
was defined in \cite{GarNhieu1996}. Let us define  a
$d_{CC}-$annulus $A^{R}_{r}:=\{x\in {\mathbb{R}}^{n}:\
r<d_{CC}(x)<R\}.$

Let us introduce ${\mathfrak L}_{1,p}(\Omega)$ for $1\leq p<\infty$
as follows: it is a weak Sobolev space, which is a norm closure of a
set of functions  $C^\infty (\bar{\Omega})$ with the norm
\begin{equation*}
  \|f\|_{{\mathfrak
L}^{1,p}(\Omega)}=\left(
 \underset{\Omega}\int \left( |\mathbf{X} f|^p + |f|^p \right)dx\right)^{\frac{1}{p}},
\end{equation*}
which is an equivalent to the norm
\begin{equation*}
  \left(
 \underset{\Omega}\int  |\mathbf{X} f|^p
 dx\right)^{\frac{1}{p}}+\left(
 \underset{\Omega}\int  |f|^r
 dx\right)^{\frac{1}{r}},
\end{equation*} where $0<r\leq p.$

${\overset{\circ}{\mathfrak L}}_{1,p}(\Omega)$ is a subspace of
${\mathfrak L}_{1,p}(\Omega)$, which is a norm closure of functions
from ${C_0^\infty}({\Omega})$ with the norm of ${\mathfrak
L}_{1,p}(\Omega)$.


It is well known that for  ${\mathfrak L}_{1,p}(\Omega)$ and
${\overset{\circ}{\mathfrak L}}_{1,p}(\Omega)$ under
Carnot-Carath\'{e}odory spaces were proved embedding Sobolev type
theorems  and Nirenberg-Galiardo type inequalities.
 We would refer the reader  interested in the theory of Carnot-Carath\'{e}odory spaces to the survey \cite{GarNhieu1996},
 where one can find more useful references. Particularly,  from the results of the survey \cite{GarNhieu1996} one can easy prove the following multiplicative Nirenberg-Galiardo type
 inequality.
 \begin{prop}\label{markp3.3}
For every function $ f \in {{\mathfrak
L}}_{1,p}({\mathbb{R}}^{N+M}), $ the following inequality holds
\begin{equation}\label{mark4.5}
   \int_{{\mathbb{R}}^{N+M}} |f|^{p} dz  \leq
  C_{14} {\left(\int_{{\mathbb{R}}^{N+M}} |\mathbf{X} f|^{p}dz \right)}^{\beta_1}
  {\left(\int_{{\mathbb{R}}^{N+M}} |f|^{\beta_2}dz \right)}^{\frac{(1-\beta_1)p}{\beta_2}}.
\end{equation}
Here $0<\beta_1=\beta_1(\beta_2)<1$.
\end{prop}

\begin{defn}
We can say that a function $u(x)$ as a generalized entire solution
(a weak solution) of the equation \eqref{mark4.1} if
$u(x)\in{\mathfrak L}_{1,p,loc}({\mathbb{R}}^{n})$ and satisfies the
inequality
\begin{equation}\label{mark4.2}
 \underset{\Omega} {\int} {\overset{m}{\underset{i=1}{\sum}}}  A_i(x,u,\mathbf{X}
 u)X_i \varphi
 +
f(x,u)\varphi  dx\leq 0,
\end{equation}for any $\varphi\in{\overset{\circ}{\mathfrak L}}_{1,p}(\Omega)$ and
$\Omega$ is an open bounded domain from ${\mathbb{R}}^{n}$.
\end{defn}


\section{Proof of a priori estimates for the case $q\leq p$}

 In this section we will prove  results
indicated by i) in Theorem \ref{markt4.1} and Theorem
\ref{markt4.2}. At first let us mention that using the simple change
of notation for $u$ it is easy to check that conditions
$(\bf{A_1})$, $(\bf{A_2})$, $(\bf{A_3})$ could be rewritten in the
following form
$$
\bf{(A_1')}\begin{array}{l} A(x,u-\tilde{b}, \xi)\  \mbox{is a
Carnot function},\\\mbox{such that for all} \ \xi\in
{\mathbb{R}}^{m}: \nu_1 |\xi|^p \leq A (x,u-\tilde{b},{\xi}){\xi}\
\mbox{on}\ {\mathbb{R}}^{n},
 \end{array}
$$
$$
\bf{(A_2')}\begin{array}{l} A(x,u-\tilde{b},\vec{\xi})\ \mbox{is a
Carnot function},\\\mbox{such that for all}\ \xi\in
{\mathbb{R}}^{m}: A (x,u-\tilde{b},{\xi})\leq \nu_2 |{\xi}|^{p-1}\
\mbox{on}\ {\mathbb{R}}^{n},
 \end{array}
$$
$$
\bf{(A_3')}\begin{array}{l} f(x,u)\ \mbox{is a Carnot
function},\\\mbox{such that for all}\ u \in {\mathbb{R}}:
f(x,u-\tilde{b}) \geq \nu_3|u-\tilde{b}|^{q-2}(u-\tilde{b})\
\mbox{on}\ {\mathbb{R}}^{n}
 \end{array}
$$where $\tilde{b}$ is a fixed positive constant.
 Let us define a cutoff
function $\xi_R$ from the space of $d_{CC}$-Lipschitz continuous
functions such that $\xi_R\equiv 0$ if $d_{CC}(x)\leq R$ or
$d_{CC}(x)\geq \mu R$ with $|\mathbf{X} \xi_R|\leq C/(\mu-1)R$ a.e.
in ${\mathbb{R}}^{n}$ for any $R>0, \mu>1.\ $ This function exists
according to \cite{GarNhieu1998}.

 Then we can use
$(u-\tilde{b})_+^s \xi_R^m$ ($\tilde{b}>0,\ m,s>1$ we shall choose
in what follows) as a test function for integral inequality
\eqref{mark4.2} which together with standard Young's inequality
yields basic for the case $p\neq q$ a priori estimate

\begin{lem}\label{markl4.1} Let $u$ be a weak solution of the
inequality  \eqref{mark4.1}.  Then the following a priori estimate
holds locally on any annulus $A_{R}^{\mu R},$ where $R>0$ is big
enough and $\mu>1$ and for any level number $\tilde{b}>0$ and
exponents $m,s>1$
\begin{equation*}
\underset{A_{R}^{\mu R}\cap\{u\geq \tilde{b} \}}\int
 (u-\tilde{b})_+^{q+s-1} \xi_R^m dx
  +\underset{A_{R}^{\mu R}\cap\{u\geq \tilde{b} \}}\int
  |\mathbf{X} (u-\tilde{b})_+^{\frac{p+s-1}{p}} |^{p} \xi_R^m dx
 \end{equation*}
\begin{equation}\label{mark4.6}
    \leq C \underset{A_{R}^{\mu R}\cap\{u\geq \tilde{b} \}}\int
    (u-\tilde{b})_+^{p+s-1}|\mathbf{X} \xi_R|^p \xi_R^{m-p} dx,
\end{equation} where $C=C(\varepsilon, m,p,s),\ m>p$ and $\varepsilon
>0$ is small enough.
\end{lem}
To investigate the case $p=q$ we will need the following
\begin{lem}\label{markl4.2}
Let $u$ be a weak solution of the inequality  \eqref{mark4.1} for
$1<q=p<Q$. Then the following a priori estimate holds locally on any
annulus $A_{R}^{\mu R},$ where  $R>0$ is big enough and $\mu>1$ and
for any level number $\tilde{b}>0$
\begin{equation*}
\underset{A_{R}^{\mu R}\cap\{u\geq \tilde{b} \}}\int
 (u-\tilde{b})_+^{p-1} \xi_R^m dx
  +\underset{A_{R}^{\mu R}\cap\{u\geq \tilde{b} \}}\int
  |\mathbf{X} (u-\tilde{b})_+^{\frac{2(p-1)}{p}} |^{p} \xi_R^m dx
 \end{equation*}
\begin{equation}\label{mark4.7}
    \leq C \underset{A_{R}^{\mu R}\cap\{u\geq \tilde{b} \}}\int
    (u-\tilde{b})_+^{2(p-1)}|\mathbf{X} \xi_R|^p \xi_R^{m-p} dx,
\end{equation} where $C=C(\varepsilon_1, m,p),\ m>p$ and
$\varepsilon_1
>0$ is small enough.
\end{lem}

Indeed, let us take
$\frac{(u-\tilde{b})_+}{(u-\tilde{b})_++\bar{\varepsilon}}\xi_R^{m}$
$(\bar{\varepsilon}>0)$ as a test function, then using the same
argument as in the previous lemma we obtain
\begin{equation*}
\underset{A_{R}^{\mu R}\cap\{u\geq \tilde{b} \}}\int \frac{
(u-\tilde{b})_+^{p} }{(u-\tilde{b})_++\bar{\varepsilon}}\xi_R^m dx
  +(\bar{\varepsilon}-\varepsilon_1)\underset{A_{R}^{\mu R}\cap\{u\geq \tilde{b} \}}\int
 \frac{ |\mathbf{X} (u-\tilde{b})_+ |^{p}}{\left((u-\tilde{b})_++\bar{\varepsilon}\right)^2} \xi_R^m dx
 \end{equation*}
\begin{equation}\label{mark4.8}
    \leq C(\varepsilon_1)m^p \underset{A_{R}^{\mu R}\cap\{u\geq \tilde{b} \}}\int
    (u-\tilde{b})_+^{p}((u-\tilde{b})_++\bar{\varepsilon})^{p-2}|\mathbf{X} \xi_R|^p \xi_R^{m-p}
    dx.
\end{equation}
Throwing away the second term from the left and letting
$\bar{\varepsilon} \rightarrow 0$ we obtain
\begin{equation} \label{mark4.9}
\underset{A_{R}^{\mu R}\cap\{u\geq \tilde{b} \}}\int {
(u-\tilde{b})_+^{p-1} }\xi_R^m dx
    \leq C(\varepsilon_1)m^p \underset{A_{R}^{\mu R}\cap\{u\geq \tilde{b} \}}\int
    (u-\tilde{b})_+^{2(p-1)}|\mathbf{X} \xi_R|^p \xi_R^{m-p}
    dx.
\end{equation}
Choosing $s=p-1$ in \eqref{mark4.6} and taking into account that
$q=p$ we can compose the required inequality from \eqref{mark4.6}
and \eqref{mark4.9}.

The assumptions of the next lemma are satisfied both for the case
$p\neq q$ and $p=q$.
\begin{lem}\label{markl4.3}
Suppose that for all weak solutions of the inequality
\eqref{mark4.1} for some fixed $l>0$ and $0<\delta<1$ the following
a priori estimate holds locally on any annulus $A_{\tilde{R}}^{\mu
\tilde{R}},$ where $\tilde{R}>0$ is big enough and $\mu>1$ and for
any level number $\tilde{b}>0$
\begin{equation*}
\underset{A_{\tilde{R}}^{\mu \tilde{R}}\cap\{u\geq \tilde{b} \}}\int
 (u-\tilde{b})_+^{l} \xi_{\tilde{R}}^m dx
  +\underset{A_{\tilde{R}}^{\mu \tilde{R}}\cap\{u\geq \tilde{b} \}}\int
  |\mathbf{X} (u-\tilde{b})_+^{\frac{l(1+\delta)}{p}} |^{p} \xi_{\tilde{R}}^m dx
 \end{equation*}
\begin{equation}\label{mark4.11}
    \leq C \underset{A_{\tilde{R}}^{\mu \tilde{R}}\cap\{u\geq \tilde{b} \}}\int
    (u-\tilde{b})_+^{l(1+\delta)}|\mathbf{X} \xi_{\tilde{R}}|^p \xi_{\tilde{R}}^{m-p}
    dx.
\end{equation}
Here $p$ is the parameter of the inequality \eqref{mark4.1} and
$p<Q$, $m$ is any fixed number greater than $2p,$
$C=C(m,p,\delta,Q)$ is some constant usually larger than $1$. Then
for any $R>0$, for any fixed $\sigma_1,\ \sigma_2,\ \sigma_3>0$ and
for any level numbers $b_1>\tilde{b}>b_2>0$ and $0<\nu<1$ one can
find proper $\beta_1\in(0,1)$ which depends only from $p,l, \delta$
and $\nu$ such that
\begin{equation}\label{mark4.12}
\underset{A_{R+\sigma_1 R}^{R+\sigma_1 R+\sigma_2 R} }\int
 (u-{b}_1)_+^{l} \xi_R^m dx \leq C
 (\sigma_{min}R)^{-\frac{p}{1-\beta_1}}\left(\underset{A_{R}^{R+\sigma_1 R+\sigma_2 R+\sigma_3 R} }\int
 (u-{b}_2)_+^{l\nu} \xi_R^m dx\right)^{\frac{1+\delta}{\nu}}
\end{equation}
\end{lem}
Suppose that $r^{l}_j=R(1+\sigma_1 2^{-j}),$ $\
r^{r}_j=R(1+\sigma_1+ \sigma_2+\sigma_3-\sigma_3 2^{-j}),$ $\
b_j=b_2+(b_1-b_2)2^{-j},$ $\ \tilde{b}_{j+1}=(b_j+b_{j+1})/2,$ $\
j=0,1,2,..$ and let $\xi_{j+1}(x)$ be a $d_{CC}$- Lipschitz
continuous function such that $\xi_{j+1}(x) \in
{\overset{\circ}{\mathfrak L}}_{1,p}({\mathbb{R}}^{n}),\
\xi_{j+1}\equiv 1 $ for all $x\in A_{r^l_j}^{r^r_j},$
$\xi_{j+1}\equiv 0$ for all $x$ if $r^r_{j+1}\leq d_{CC}(x)$ or when
$d_{CC}(x)\leq r^l_{j+1},$ $|\mathbf{X} \xi_{j+1}|=2^{j}|\mathbf{X}
\xi_{j} |, $ $|\mathbf{X} \xi_{j+1}|\leq C/(r^l_j-r^l_{j+1})$ for a.
e.  $x \in {\mathbb{R}}^{n}.$ This function exists according to
\cite{GarNhieu1998}. Setting in the inequality \eqref{mark4.11}
$\tilde{R}=r^l_{j+1}, \mu \tilde{R}=r^r_{j+1},
\tilde{b}=\tilde{b}_{j+1},  \xi_{\tilde{R}}=\xi_{j+1},
A_{r^l_{j+1}}^{r^r_{j+1}}=A_{j+1}, $ we obtain

\begin{equation*}
\underset{A_{j+1}\cap\{u\geq {\tilde{b}}_{j+1} \}}\int
 (u-{\tilde{b}}_{j+1})_+^{l} \xi_{j+1}^m dx
  +\underset{A_{j+1}\cap\{u\geq {\tilde{b}}_{j+1}}\int
  |\mathbf{X} (u-{\tilde{b}}_{j+1})_+^{\frac{l(1+\delta)}{p}} |^{p} \xi_{j+1}^m dx
 \end{equation*}
\begin{equation}\label{mark4.13}
    \leq C \underset{A_{j+1}\cap\{u\geq {\tilde{b}}_{j+1} \}}\int
    (u-{\tilde{b}}_{j+1})_+^{l(1+\delta)}|\mathbf{X} \xi_{j+1}|^p \xi_{j+1}^{m-p}
    dx.
\end{equation}
Denote $f_j=(u-b_j)_+^{\frac{l(1+\delta)}{p}}
\xi_j^{\frac{m-p}{p}}.$ Then using this and previous notations and
listed above properties of cut-off functions $\xi_j$ and
$\xi_{j+1}$, we can get
\begin{equation*}
\underset{A_{j}}\int
 |f_j|^{\frac{p}{1+\delta}}  dx
  +\underset{A_{j}}\int
  |\mathbf{X}f_j |^{p}  dx \leq \underset{A_{j+1}\cap\{u\geq {\tilde{b}}_{j+1} \}}\int
 (u-{\tilde{b}}_{j+1})_+^{l} \xi_{j+1}^m dx
\end{equation*}
\begin{equation}\label{mark4.14}
+\underset{A_{j+1}\cap\{u\geq {\tilde{b}}_{j+1}}\int
  |\mathbf{X} (u-{\tilde{b}}_{j+1})_+^{\frac{l(1+\delta)}{p}} |^{p} \xi_{j+1}^m
  dx+ \end{equation}
\begin{equation*}\left(\frac{m-p}{p}\right)^p \underset{A_{j}}\int
    (u-{\tilde{b}}_{j})_+^{l(1+\delta)}|\mathbf{X} \xi_{j+1}|^p \xi_{j+1}^{m-p}
    dx.
\end{equation*}
Then it follows from \eqref{mark4.13} and \eqref{mark4.14} that
\begin{equation}\label{mark4.15}
\underset{A_{j}}\int
 |f_j|^{\frac{p}{1+\delta}}  dx
  +\underset{A_{j}}\int
  |\mathbf{X}f_j |^{p}  dx \leq  C  \frac{2^{pj}}{(\sigma_{min}R)^p}\underset{A_{j+1}}\int
  |f_{j+1}|^p dx
\end{equation}
To estimate from above the integral on the right hand side of the
inequality \eqref{mark4.15} we need to use the Nirenberg-Gagliardo
type multiplicative inequality \eqref{mark4.5}. Let us apply it to
the function $f=f_{j+1}$ with $\beta_2=p/(1+\delta).$ This yields
\begin{equation}\label{mark4.16}
\underset{A_{j}}\int
 |f_j|^{\frac{p}{1+\delta}}  dx
  +\underset{A_{j}}\int
  |\mathbf{X}f_j |^{p}  dx \leq
\end{equation}
\begin{equation*}
   C  \frac{2^{pj}}{(\sigma_{min}R)^p}
  \left(\underset{A_{j+1}}\int
  |\mathbf{X}f_{j+1} |^{p}  dx \right)^{\beta_1}\left( \underset{A_{j+1}}\int
  |f_{j+1}|^{\frac{p\nu}{1+\delta}}
  dx\right)^{\frac{(1-\beta_1)(1+\delta)}{\nu}},
\end{equation*}
where, after calculations,
$\beta_1=\frac{\frac{1+\delta}{\nu}-1}{\frac{1+\delta}{\nu}-1+\frac{p}{Q}}<1.$

Using the standard Young's inequality with exponents $\beta_1^{-1}$
and $(1-\beta_1)^{-1}$ together with \eqref{mark4.15}, we have
\begin{equation*}
\underset{A_{j}}\int
 |f_j|^{\frac{p}{1+\delta}}  dx
  +\underset{A_{j}}\int
  |\mathbf{X}f_j |^{p}  dx \leq \varepsilon_2   \underset{A_{j+1}}\int
  |\mathbf{X}f_{j+1} |^{p}  dx + \end{equation*}
\begin{equation*} \frac{ C(\varepsilon_2) \left(2^{\frac{p}{1-\beta_1}}\right)^j}{(\sigma_{min}R)^{\frac{p}{1-\beta_1}}}
\left( \underset{A_{j+1}}\int
  |f_{j+1}|^{\frac{p\nu}{1+\delta}}
  dx\right)^{\frac{1+\delta}{\nu}},
\end{equation*}where $\varepsilon_2>0$ is a small constant to be
chosen later. By induction we obtain
\begin{equation*}
\underset{A_{0}}\int
 |f_0|^{\frac{p}{1+\delta}}  dx
  +\underset{A_{0}}\int
  |\mathbf{X}f_0 |^{p}  dx \leq \varepsilon_2^{j+1}   \underset{A_{j+1}}\int
  |\mathbf{X}f_{j+1} |^{p}  dx +
  \end{equation*}
\begin{equation*} \frac{ C(\varepsilon_2)\sum_{k=0}^{j}\left(\varepsilon_2 2^{\frac{p}{1-\beta_1}}\right)^k
  }{(\sigma_{min}R)^{\frac{p}{1-\beta_1}}}
\left( \underset{A_{j+1}}\int
  |f_{j+1}|^{\frac{p\nu}{1+\delta}}
  dx\right)^{\frac{1+\delta}{\nu}}.
\end{equation*}
Here  we can cast aside the gradient term from the left. Then
according to the definition of weak solutions the integral
$\int_{A_{j+1}}
  |\mathbf{X}f_{j+1} |^{p}  dx$ is convergent. Choosing
  $\varepsilon_2$ so that $\varepsilon_2 2^{\frac{p}{1-\beta_1}}=1/2<1$
  and letting $j\rightarrow\infty,$ we obtain the required estimate.

  \begin{lem}\label{markl4.4}
  Under the assumptions of the previous lemma the following
  inequality holds
  \begin{equation}\label{mark4.17}
  \|u\|_{\infty,\
  A_{\bar{R}+\bar{\sigma}_1\bar{R}}^{\bar{R}+\bar{\sigma}_1\bar{R}+\bar{\sigma}_2\bar{R}}}\leq
  C ({\bar{\sigma}}_{min}\bar{R})^{-\frac{p}{1-\beta_1}\cdot
  \frac{1}{l(1-\nu)}} \left(\underset{A_{\bar{R}}^{\bar{R}+{\bar{\sigma}}_1 \bar{R}+{\bar{\sigma}}_2 \bar{R}+{\bar{\sigma}}_3 \bar{R}} }\int
 |u|^{l\nu}
 dx\right)^{\frac{1}{l\nu}\cdot\frac{1+\delta-\nu}{1-\nu}}.
\end{equation}
\end{lem}
Suppose that $R_i^l=\bar{R}(1+\bar{\sigma}_1-\bar{\sigma}_1
2^{-i}),\ {\tilde{R}}_i^l=(R_i^l+R_{i+1}^l)/2,\
R_i^r=\bar{R}(1+\bar{\sigma}_1+\bar{\sigma}_2+\bar{\sigma}_3
2^{-i}),\ {\tilde{R}}_i^r=(R_i^r+R_{i+1}^r)/2,\
{\bar{\sigma}}_{min}=\min(\bar{\sigma}_1;\bar{\sigma}_3),\
h_i=k(1-2^{-i-1}),\ {\tilde{h}}_i=(h_i+h_{i+1})/2\ $ for $\
i=0,1,2..$ Then under the settings $R=R_i^l, \ R+\sigma_1
R={\tilde{R}}_i^l, \ R+\sigma_1 R+\sigma_2 R={\tilde{R}}_i^r, \
R+\sigma_1 R+\sigma_2 R+\sigma_3 R=R_i^r, \ b_1={\tilde{h}}_i,\
b_2=h_i, \ A_{R+\sigma_1 R}^{R+\sigma_1 R+\sigma_2
R}=A_{{\tilde{R}}_i^l}^{{\tilde{R}}_i^r}={\tilde{A}}_i,
A_R^{R+\sigma_1 R+\sigma_2 R+\sigma_3 R}=A_{R_i^l}^{R_i^r}=A_i$
Lemma \ref{markl4.3} implies
\begin{equation}\label{mark4.18}
\underset{{\tilde{A}}_i\cap\{u\geq {\tilde{h}}_i \}}\int
 (u-{\tilde{h}}_i)_+^{l}  dx \leq C ({\tilde{R}}_i^l-R_i^l)^{-\frac{p}{1-\beta_1}}\left(\underset{A_i\cap\{u\geq {\tilde{h}}_i \} }\int
 (u-h_i)_+^{l\nu}  dx\right)^{\frac{1+\delta}{\nu}}
\end{equation}
Let us denote $I_{i+1}=\underset{{{A}}_{i+1}}\int
 (u-{{h}}_{i+1})_+^{l\nu}  dx.$ Thus, from \eqref{mark4.18} one can
 get the following estimate
 \begin{equation*}
 I_{i+1}\leq C ({\tilde{R}}_i^l-R_i^l)^{-\frac{p}{1-\beta_1}}
 (h_{i+1}-h_i)^{-l(1-\nu)}I_i^{\frac{1+\delta}{\nu}}.
\end{equation*}
It means that the sequence $\{I_i\}$ satisfies the assumptions of
Ladizhenskaya's Lemma and from the last inequality we have
$I_{i+1}\leq C_{Lad} b^i I_i^{1+\theta},$ denoting
$\theta={\frac{1+\delta}{\nu}}-1>0,$
$$
C_{Lad}=C
({\bar{\sigma}}_{min}\bar{R})^{-\frac{p}{1-\beta_1}}(2^{{\frac{p}{1-\beta_1}}+l(1-\nu)})^{i}k^{-l(1-\nu)}>0,\
b=2^{{\frac{p}{1-\beta_1}}+l(1-\nu)}>1,$$ where the constant
$C_{Lad}$ is controlled by $k$. Then Ladizhenskay's Lemma implies
that $I_i\rightarrow 0$ as $i\rightarrow\infty$, and $\|u\|_{\infty,
A_{\infty}}\leq k$ if $I_0\leq C_{Lad}^{-1/\theta}b^{-1/\theta^2}.$
Choosing $$k=C
({\bar{\sigma}}_{min}\bar{R})^{-\frac{p}{1-\beta_1}\cdot
  \frac{1}{l(1-\nu)}} \left(\underset{A_{\bar{R}}^{\bar{R}+{\bar{\sigma}}_1 \bar{R}+{\bar{\sigma}}_2 \bar{R}+{\bar{\sigma}}_3 \bar{R}} }\int
 |u|^{l\nu}
 dx\right)^{\frac{1}{l\nu}\cdot\frac{1+\delta-\nu}{1-\nu}},$$ where
 constant $C$ is sufficiently large, we obtain the required
 estimate. The lemma \ref{markl4.4} is proved.

To make a simple iteration process, redenoting manifest that
together with Lemma \ref{markl4.4}, we have the following inequality
on any sequence  of nested annuli $A_i$ such that $A_{\mu R}^{2\mu
R}=A_\infty\subset A_{i+1}\subset  A_{i}\subset A_0=A_R^{3\mu R}$
\begin{equation}\label{mark4.19}
\|u\|_{\infty, A_{i+1}}\leq C b^i R^{-\frac{p}{1-\beta_1}\cdot
  \frac{1}{l(1-\nu)}} \left(\underset{A_i }\int
 |u|^{l\nu}
 dx\right)^{\frac{1}{l\nu}\cdot\frac{1+\delta-\nu}{1-\nu}}\leq \end{equation}
\begin{equation*} C b^i R^{-\frac{p}{1-\beta_1}\cdot
  \frac{1}{l(1-\nu)}}R^{\frac{n(\frac{1+\delta}{\nu}-1)}{l(1-\nu)}} \|u\|_{\infty, A_i}^{1+
  \frac{\delta}{1-\nu}}.
  \end{equation*}
Let us denote also $c=C R^{-\frac{p}{1-\beta_1}\cdot
  \frac{1}{l(1-\nu)}}R^{\frac{n(\frac{1+\delta}{\nu}-1)}{l(1-\nu)}}, \ \varepsilon=\frac{\delta}{1-\nu}.
  $ Thus,
  $$
\|u\|_{\infty, A_{i+1}}\leq c b^i\|u\|_{\infty,
A_i}^{1+\varepsilon},
  $$
$$
\|u\|_{\infty, A_{\infty}}\leq\|u\|_{\infty, A_{i}}\leq
c^{\frac{(1+\varepsilon)^i-1}{\varepsilon}}b^{\frac{(1+\varepsilon)^i-1}{\varepsilon^2}-\frac{i}{\varepsilon}}\|u\|_{\infty,
A_0}^{(1+\varepsilon)^i},
$$
$$
\|u\|_{\infty,A_{\mu R}^{2\mu R}}^{\frac{1}{(1+\varepsilon)^i}}\leq
c^{\frac{1}{\varepsilon}-\frac{1}{\varepsilon(1+\varepsilon)^i }}
b^{\frac{1}{\varepsilon^2}-\frac{1}{\varepsilon^2(1+\varepsilon)^i
}-\frac{i}{\varepsilon(1+\varepsilon)^i}}\|u\|_{\infty, A_R^{3\mu
R}}.
$$
As far as $\|u\|_{\infty,A_{\mu R}^{2\mu R}}\neq 0$ is a fixed
number for any fixed $R$,   we obtain as $i\rightarrow \infty$

$$
1\leq c^{\frac{1}{\varepsilon}}
b^{\frac{1}{\varepsilon^2}}\|u\|_{\infty, A_R^{3\mu R}},
$$
$$
1\leq c b^{\frac{1}{\varepsilon}}\|u\|_{\infty, A_R^{3\mu
R}}^\varepsilon.
$$

 This yields
\begin{equation}\label{mark4.21}
 R^{\frac{p}{1-\beta_1}\cdot
\frac{1}{l(1-\nu)}-\frac{n(\frac{1+\delta}{\nu}-1)}{l(1-\nu)}} \leq
C(\mu, n)\|u\|_{\infty, A_R^{3\mu R}}^{\frac{\delta}{1-\nu}},
\end{equation}
and after the appropriate calculations, from \eqref{mark4.21} we
have the desired result
$$
R^{\frac{p}{l\delta}}\leq C \|u\|_{\infty, A_R^{3\mu R}}.
$$
Now we need only to choose appropriate parameters such as $l$ and
$\delta$ according to $q<p$ case and for the case $q=p$.

\section{Proof of Theorem \ref{markt4.3}}

Here we will get the a priori estimate with $p<q.$ Let $u$ be a weak
solution of the inequality \eqref{mark4.1}. Following the same way
as it was in the proof of Lemma \ref{markl4.1} and using
non-negativity of characteristic form of our differential operator,
we can obtain
\begin{equation*}
\underset{\Omega}\int
 |u|^{q+s-1} \xi_R^m dx
\label{mark4.50}
    \leq \frac{ m^p}{s 2^p }\underset{\Omega}\int
    |u|^{p+s-1}|\mathbf{X} \xi_R|^p \xi_R^{m-p} dx.
\end{equation*}It gives  us together with Holder inequality the
following estimate
\begin{equation*}
\underset{\Omega}\int
 |u|^{q+s-1} \xi_R^m dx
\label{mark4.51}
    \leq \frac{ m^p}{s 2^p (R-r)^p}C(Q)R^{\frac{Q(q-p)}{q+s-1}}\left(\underset{\Omega}\int
    |u|^{q+s-1} \xi_R^{\frac{(m-p)(q+s-1)}{p+s-1}}
    dx\right)^{\frac{p+s-1}{q+s-1}}.
\end{equation*}Let us fix $m=\frac{p(q+s-1)}{q-p}$ such that it satisfies
$m={\frac{(m-p)(q+s-1)}{p+s-1}}$ then integrals from both sides are
equal. Hence, we have
\begin{equation*}
\left(\underset{B_R}\int
    |u|^{q+s-1} \xi_R^{\frac{p(q+s-1)}{q-p}}
    dx\right)^{\frac{q-p}{q+s-1}}\leq \frac{m^p C(Q)R^{\frac{Q(q-p)}{q+s-1}}}{s 2^p
    (R-r)^p}=\end{equation*}
\begin{equation*} C
    \frac{R^{\frac{Q(q-p)}{q+s-1}}}{(R-r)^p}\frac{p^p}{(q-p)^p}(1+\frac{q-1}{s})(q+s-1)^{p-1}.
    \end{equation*}
This yields
\begin{equation*}
\left(\underset{B_R}\int
    |u|^{q+s-1} \xi_R^{\frac{p(q+s-1)}{q-p}}
    dx\right)^{\frac{q-p}{q+s-1}}\leq \frac{C(p,q)R^{\frac{Q(q-p)}{q+s-1}} (q+s-1)^{p-1}}{(R-r)^p}(1+\frac{q-1}{s}).
    \end{equation*}
Let us choose $r=R/2$ and redenote $\bar{s}=q+s-1$. Then we have
\begin{equation*}
\left(\underset{B_{R/2}}\int
    |u|^{\bar{s}}
    dx\right)^{\frac{1}{\bar{s}}}\leq C(p,q)
\left(\frac{R^{\frac{Q(q-p)}{\bar{s}}}{\bar{s}}^{{p-1}}}{R^{{p}}}\right)^{\frac{1}{q-p}}
(1+\frac{q-1}{\bar{s}-q+1})^{\frac{1}{q-p}} .
    \end{equation*}Now, as $R\rightarrow +\infty$ then for all $\bar{s}> \frac{Q(q-p)}{p}$ we have $L_{\bar{s}}$-norms
    on ${\mathbb{R}}^{n}$ are zero.  Thus, Theorem \ref{markt4.3} is proved.


\subsection*{Acknowledgment}
Many thanks to professor Giovanna Citti for useful discussions.

\end{document}